\documentclass[11pt,bezier]{article}
\usepackage{amsmath}
\usepackage{amsfonts,amsthm,amssymb}
\usepackage{amsfonts}
\usepackage{graphics}
\usepackage{amssymb}
\newtheorem{lem}{Lemma}[section]
\newtheorem{thm}[lem]{Theorem}

\newtheorem{conj}{Conjecture}

\theoremstyle{definition}

\begin{document}
\title{Connectivity keeping paths in $k$-connected bipartite graphs\footnote{The research is supported by NSFC (11861066), Tianshan Youth Project of Xinjiang (2018Q066).}}
\author{Lian Luo, Yingzhi Tian\footnote{Corresponding author. E-mail: 2583224265@qq.com (L. Luo), tianyzhxj@163.com (Y. Tian), wly95954@163.com (L. Wu).}, Liyun Wu \\
{\small College of Mathematics and System Sciences, Xinjiang
University, Urumqi, Xinjiang 830046, PR China}}

\date{}

\maketitle

\noindent{\bf Abstract } In 2010,  Mader [W. Mader, Connectivity keeping paths in $k$-connected graphs, J. Graph Theory 65 (2010) 61-69.] proved that  every $k$-connected graph $G$ with minimum degree at least $\lfloor\frac{3k}{2}\rfloor+m-1$ contains a path $P$ of order $m$ such that $G-V(P)$ is still $k$-connected. In this paper, we consider similar problem for bipartite graphs, and prove that every $k$-connected bipartite graph $G$ with minimum degree at least $k+m$ contains a path $P$ of order $m$ such that $G-V(P)$ is still $k$-connected.

\noindent{\bf Keywords:} Connectivity; Bipartite graphs; Paths

\section{Introduction}

For graph-theoretical terminology and notations not defined here, we follow \cite{Bondy}. In this paper, we consider finite simple graphs.  Let $G$ be a graph with vertex set $V(G)$ and edge set $E(G)$. The $order$ of $G$ is the number of its vertices, denoted by $|G|$.  The $minimum$ $degree$ of $G$ is denoted by $\delta(G)$. The $connectivity$ of $G$, denoted by $\kappa(G)$, is the minimum size of a  vertex set $S$  such that $G-S$ is disconnected or has only one vertex. We say $G$ is $k$-connected if $\kappa(G)\geq k$. The $floor$ of a real number $x$, denoted by $\lfloor x\rfloor$, is the greatest integer not larger than $x$; the $ceiling$ of a real number $x$, denoted by $\lceil x\rceil$, is the least integer greater than or equal to $x$.

In 1972, Chartrand, Kaugars, and Lick proved the following well-known result.

\begin{thm}{\cite{Chartrand}} Every $k$-connected graph $G$ with minimum degree
$\delta(G)\geq \lfloor\frac{3}{2}k\rfloor$ has a vertex $v$ such that $G-v$ is still $k$-connected.
\end{thm}

In 2003, Fujita and Kawarabayashi \cite{Fujita} considered a similar problem for an edge of a graph, and proved that every $k$-connected graph $G$ with minimum degree at least $\lfloor\frac{3}{2}k\rfloor+2$ has an edge $e=uv$ such that $G-\{u,v\}$ is still $k$-connected. They also proposed the following conjecture.

\begin{conj}{\cite{Fujita}} For all positive integers $k, m$, there
is a (least) non-negative integer $f_k(m)$ such that every $k$-connected graph $G$ with $\delta(G)\geq \lfloor\frac{3}{2}k\rfloor-1+f_k(m)$ contains a connected subgraph $W$ of exact order $m$ such that $G-V(W)$ is still $k$-connected.
\end{conj}

The examples given in \cite{Fujita} showed that $f_k(m)$ must be at least $m$ for all positive integers $k, m$. In \cite{Mader1}, Mader confirmed Conjecture 1 and  proved that $f_k(m)=m$ holds for all $k, m$. Moreover, the connected subgraph $W$ could even be a path.

\begin{thm}{\cite{Mader1}}
Every $k$-connected graph $G$ with $\delta(G)\geq\lfloor\frac{3}{2}k\rfloor+m-1$
for positive integers $k, m$ contains a path $P$ of order $m$ such that $G-V(P)$ remains $k$-connected.
\end{thm}

Mader \cite{Mader1} further conjectured that Theorem 1.2 holds for all trees with order $m$.

\begin{conj}{\cite{Mader1}} For every positive integer $k$ and every finite tree $T$ with order $m$,  every $k$-connected graph $G$ with $\delta(G)\geq \lfloor\frac{3}{2}k\rfloor+m-1$ contains a subgraph $T'\cong T$ such that $G-V(T')$ is still $k$-connected.
\end{conj}

In \cite{Mader2}, Mader showed that Conjecture 2 is true  if $\delta(G)\geq2(k-1+m)^2+m-1$. The result in \cite{Diwan} showed that Conjecture 2 is true for $k=1$. For general $k\geq2$, Conjecture 2 remains open. But for $k=2$, partially affirmative answers have been shown, see [6-9,12,13].

Motivated by Theorem 1.2, we consider a similar problem for bipartite graphs in this paper. In the next section, we will present notations, terminology and some lemmas used in the proofs of main results.   Section 3 gives main results of this paper.  Some remarks will be concluded in the last section.

\section{Preliminaries}

For each vertex $v\in V(G)$, the $neighborhood$ $N_G(v)$ of $v$ is defined as the set of all vertices adjacent to $v$,  and $d_G(v)=|N_G(v)|$ is the $degree$ of $v$. For a subgraph $H\subseteq G$, we define $\delta_G(H)=$min$_{v\in V(H)}d_G(v)$, whereas $\delta(H)$ is the minimum degree of the graph $H$. For $S\subseteq V(G)$, we set $N_G(S)=(\cup_{v\in S}N_G(v))\setminus S$. The induced subgraph of $S$ in $G$ is denoted by $G[S]$. And $G-S$ is the induced subgraph $G[V(G)\setminus S]$. For $H\subseteq G$, we let $H$ stand for the graph $H$ itself and its vertex set $V(H)$, so we use  $N_G(H)$,  $G[H]$ and  $G-H$  for $N_G(V(H))$, $G[V(H)]$ and $G-V(H)$, respectively.

For two graphs $G_1=(V_1,E_1)$ and $G_2=(V_2,E_2)$, the $union$ of $G_1$ and $G_2$ is  $G_1\cup G_2=(V_1\cup V_2,E_1\cup E_2)$. The $join$ of them is denoted by $G_1\vee G_2$ consisting of $G_1\cup G_2$ and all edges joining each vertex in $V_1$ and each vertex in $V_2$. In other words, the join of them can be obtained by connecting each vertex of $G_1$ to all vertices of $G_2$.
For $uv\in E(G)$ and $H\subseteq G$, $H\cup \{uv\}$ is the subgraph of $G$ with vertex set $V(H)\cup\{u,v\}$ and edge set $E(H)\cup\{uv\}$. An $x,y$-path $P$ is a path with endvertices $x$ and $y$. For $x=y$, consider $P$ as a path of length 0. For $u,v\in V(P)$, $P[u,v]=P[v,u]$ is the subpath of $P$ between $u$ and $v$, and define $P[u,v)=P[u,v]-v$ and $P(u,v]=P[u,v]-u$.

For $S\subseteq V(G)$, we call $S$ a $vertex$-$cut$ of $G$ if $G-S$ is disconnected. A $minimum$ $vertex$-$cut$ $S$ of $G$ is a vertex-cut with $|S|=\kappa(G)$. For a vertex-cut $S$ of $G$, the union of at least one, but not all components of $G-S$ is called a $semifragment$ to $S$, and $\bar{F}=G-(S\cup V(F))$ for a semifragment $F$ to $S$ is the $complementary$ $semifragment$ of $F$ to $S$ in $G$.  If $S$ is a minimum vertex-cut of $G$, a semifragment to $S$ is called a $fragment$ of $G$ to $S$. A $minimal$ $fragment$ of $G$ is a fragment which does not properly contain another fragment of $G$. If $G$ is not complete, then there are at least two minimal fragments. For a fragment $F$ of $G$ to $S$ and a fragment $F_1$ of $G$ to $S_1$, define $S(F,F_1)=(S\cap F_1)\cup(S\cap S_1)\cup(S_1\cap F)$. For a fragment $F$ of $G$ to $S$, we have $S=N_G(F)$, whereas for a semifragment $F$ of $G$ to $S$, we have only $N_G(F)\subseteq S$. Nevertheless, we will use the same notation $S(F,F_1)$ for a semifragment $F$ of $G$ to $S$ and a semifragment $F_1$ of $G$ to $S_1$.

For a set $S$, $K(S)$ denotes the complete graph on vertex set $S$. The $completion$ of $S\subseteq V(G)$ in $G$, denoted by $G\langle S\rangle$, is the graph $G\cup K(S)$.

\begin{lem}{[5,10]}
Let $G$ be a $k$-connected graph and let $S$ be a minimum vertex-cut  of $G$. Then for every fragment $F$ of $G$ to $S$, $G\langle S\rangle-V(\bar{F})$ is $k$-connected. Furthermore, if $F$ is a minimal fragment of $G$ with $\vert F\vert\geq 2$, then $G\langle S\rangle-V(\bar{F})$ is $(k+1)$-connected.
\end{lem}

%\begin{lem}{[5,10]}
%Let $G$ be a graph with $\kappa(G)=k$ and let $F_i$ be a fragment of $G$ to $S_i$ for $i=1,2$. Then we have the following properties.
%
%$(\alpha)$
%If $F_1$ is minimal fragment of $G$ with $\vert F_1\vert\geq 2$, then $G\langle S_1\rangle-V(\bar{F_1})$ is $(k+1)$-connected.
%
%$(\beta)$ If $F_1\cap F_2\ne\emptyset$ holds, then $\vert S(F_1,F_2)\vert\geq k$ follows.
%
%$(\gamma)$ If $\vert S(F_1,F_2)\vert\geq k$ holds, then $\vert F_1\cap S_2\vert\geq \vert\bar{F_2}\cap S_1\vert$ follows.
%\end{lem}
%\par
%\vskip 0.2cm

\begin{lem} \cite{Mader1}
Let $S$ be a vertex-cut of the graph $G$ with $\vert S\vert=k$ and let $S_1$ be a vertex-cut of the graph $G$ with $\vert S_1\vert=k-1$. Assume $F$ is a semifragment of $G$ to $S$ and $F_1$ is a semifragment of $G$ to $S_1$.
Furthermore, we assume $\kappa(G\langle S\rangle-V(F))\ge k$ and $\kappa(G\langle S\rangle-V(\bar{F}))\ge k$. Then we have the following properties.

$(\alpha)$ If $F\cap F_1\ne\emptyset$ holds, then $\vert S(F,F_1)\vert\geq k$ follows.

$(\beta)$ If $\vert S(F,F_1)\vert\geq k$ holds, then we have $\vert S_1\cap F\vert\geq \vert S\cap\bar{ F_1} \vert$, $\vert S\cap F_1\vert> \vert S_1\cap\bar{ F} \vert$ and $\bar{F}\cap\bar{F_1}=\emptyset$.
\end{lem}

We denote  a bipartite graph $G$ with bipartition ($X$,$Y$) by $G=G[X,Y]$.
A $complete$ $bipartite$ $graph$  is a  bipartite graph  where each vertex of $X$ is adjacent to each vertex of $Y$. A complete bipartite graph with partitions of size $|X|=r$  and $|Y|=s$, is denoted by $K_{r,s}$.

Let $\mathcal{K}_{k}^b(m)$ denote the class of all pairs $(G,C)$, where $G$ is a $k$-connected graph and $C$ is a complete subgraph of $G$ with $|C|=k$,  $G-C$ is a bipartite graph with $\delta_{G}{(G-C)}\ge k+m$ and every pair of adjacent vertices in $V(G-C)$ have no common neighbors in $G$. Let $\mathcal{K}_{k^+}^{b}(m)$ consist of all $(G,C)\in \mathcal{K}_{k}^b(m)$ with $\kappa {(G)}\ge k+1$.

%\begin{lem}
%Let $G$ be a bipartite graph with $\delta(G)\geq k+m$. Assume $S$ is a vertex-cut of $G$ with $|S|=k$ and $F$ is a semifragment of $G$ to $S$. Then  $|F|\geq k+2m$ holds.
%\end{lem}
%
%\noindent{\bf Proof.} Since $\delta(G)\geq k+m$ and $|S|=k$, there is an edge $uv\in E(G[F])$. Thus $|F|\geq|(N_G(u)\cup N_G(u))\setminus S|\geq k+m+k+m-k=k+2m$.

\begin{lem}
Let $G$ be a graph with $\kappa(G)=k$, $S$ be a minimum vertex-cut of $G$, and $F$ be a fragment of $G$ to $S$. Suppose that there is a path $P$ in $G-(S\cup V(F))$ of order at most $m$ satisfying $\kappa(G\langle S\rangle-V(F\cup P))\ge k$. Then $\kappa(G-V(P))\ge k$ if one of the following holds:

$(\alpha)$ $G$ is a bipartite graph with $\delta(G)\geq k+m$;

$(\beta)$ $(G,C)\in\mathcal{K}_k^b(m)$ and $C\subseteq G[F\cup S]$.
\end{lem}

\noindent{\bf Proof.} Since each pair of adjacent vertices in a bipartite graph $G$ has no common neighbors, we have  $|V(G)|\geq2(k+m)$ if $(\alpha)$ holds. For $(G,C)\in\mathcal{K}_k^b(m)$, since each pair of adjacent vertices in $V(G-C)$ has no common neighbors in $G$ and $\delta_{G}{(G-C)}\ge k+m$, we also have $|V(G)|\geq2(k+m)$.

If $k=1$, then $G\langle S\rangle=G$ and the results are obvious. So we assume $k\geq2$ in the following.

Define $G_P=G-V(P)$ and assume $G_P$ is not $k$-connected. Then $G_P$ has a vertex-cut $S_1$ with $\vert S_1\vert=k-1$ by $|V(G_P)|\geq2k+m$. Let $F_1$ be a semifragment of $G_P$ to $S_1$ and $\bar{F_1}$ be the complementary semifragment of $F_1$  to $S_1$ in $G_P$. Since each vertex $v$ in $G-V(F\cup S)$ has at least $m$ neighbors in $G-V(F\cup S)$, we obtain that $G_P-V(F\cup S)$ is not empty. Thus $S$ is a vertex-cut of $G_P$ with semifragment $F$. Define the complementary semifragment of $F$ to $S$ in $G_P$ as $\bar{F}_{G_P}=G_P-(S\cup V(F))=G-(V(F)\cup S\cup V(P))$. Since $G_P\langle S\rangle-V(F)=G\langle S\rangle-V(F\cup P)$ is $k$-connected by assumption on $P$ and $G_P\langle S\rangle-V(\bar{F}_{G_P})=G\langle S\rangle[S\cup F]$ is $k$-connected by applying Lemma 2.1 to the fragment $F$ of $G$ to $S$, we will use Lemma 2.2 to $G_P$.

If $(\alpha)$ holds, then since $G$ is bipartite and $P$ is a path, each vertex in $V(G)-V(P)$ has at most $\lceil\frac{|P|}{2}\rceil$ neighbors in $V(P)$. Thus  $\delta (G_P)\geq k+m-\lceil\frac{|P|}{2}\rceil\ge k+\lfloor\frac{m}{2}\rfloor$. As we have observed, $N_G(F)=S$. Take any vertex $u\in F$. Then $u$ has at least $m$ neighbors in $F$ since $d_{G[F]}(u)\geq k+m-\vert S\vert=m$. Take $v\in N_{G[F]}(u)$. Then $N_{G_P}(u)\cap N_{G_P}(v)=\emptyset$ since they are in different partite sets. Suppose $V(F)\subseteq S_1$. Then $N_{G_P}(u)\cup N_{G_P}(v)\subseteq S\cup V(F)\subseteq S\cup S_1$. Therefore, $d_{G_P}(u)+d_{G_P}(v)\le \vert S\vert+\vert S_1\vert=2k-1$, which contradicts to $d_{G_P}(u)+d_{G_P}(v)\ge2(k+\lfloor\frac{m}{2}\rfloor)$. Then $V(F)\nsubseteq S_1$.  Thus $F\cap F_1\ne\emptyset$ or $F\cap \bar{F_1}\ne\emptyset$. Without loss of generality, we assume $F\cap F_1\ne\emptyset$. By Lemma 2.2, $\bar{F}_{G_P}\cap\bar{F_1}=\emptyset$.

If $(\beta)$ holds, then since $G-V(C)$ is bipartite and $P$ is a path, each vertex in $V(G)-V(C)$ has at most $\lceil\frac{|P|}{2}\rceil$ neighbors in $V(P)$. Thus $\delta_{G_P} (G_P-V(C))\geq k+m-\lceil\frac{|P|}{2}\rceil\ge k+\lfloor\frac{m}{2}\rfloor$. Since $(G,C)\in\mathcal{K}_k^b(m)$, $G$ in not complete. Therefore there exists a fragment of $G$ to $S$. Thus we may guarantee the existence of the fragment $F$ of $G$ to $S$ satisfying $C\subseteq G[F\cup S]$. Then $V(\bar{F}_{G_P})\cap V(C)=\emptyset$. In addition, since $C$ is complete, we can choose $F_1$ in the third paragraph of this proof such that $V(\bar{F_1})\cap V(C)=\emptyset$. Therefore $V(C)\subseteq V((F\cap F_1)\cup S(F,F_1))$. Thus  $V(C)\cap(F\cap F_1)\neq\emptyset$ or $V(C)\subseteq S(F,F_1)$. In both cases $|S(F,F_1)|\geq k$ holds by Lemma 2.2($\alpha$) or $|C|=k$. Then Lemma 2.2($\beta$) implies $\bar{F}_{G_P}\cap\bar{F_1}=\emptyset$.

Suppose, to the contrary, that $F\cap \bar{F_1}=\emptyset$ if $(\alpha)$ holds;     $|S(F,\bar{F_1})|<k$ if $(\beta)$ holds ( by Lemma 2.2$(\alpha)$, we have $F\cap \bar{F_1}=\emptyset$). Thus $\bar{F_1}\subseteq S$. Take $u\in \bar{F_1}$, since $N_G(u)\subseteq S_1\cup V(P)\cup V(\bar{F_1})$ and $u$ has at most $\lceil\frac{|P|}{2}\rceil$ neighbors in $V(P)$, we have     $d_{G[\bar{F_1}]}(u)\geq k+m-\lceil\frac{|P|}{2}\rceil-|S_1| \geq\lfloor\frac{m}{2}\rfloor+1$. Choose an edge $uv\in E(G[\bar{F_1}])$. If $(\alpha)$ holds, then $N_{G_P}(u)\cap N_{G_P}(v)=\emptyset$ is obvious. If $(\beta)$ holds,  then since $V(C)\subseteq V((F\cap F_1)\cup S(F,F_1))\subseteq F_1\cup S_1$, we have $V(\bar{F_1})\cap V(C)=\emptyset$ and $\{u,v\}\subseteq V(G)\setminus V(C)$. Then $N_{G_P}(u)\cap N_{G_P}(v)=\emptyset$ by the definition of $\mathcal{K}_k^b(m)$. Moreover,  $N_{G_P}(u)\cup N_{G_P}(v)\subseteq S_1\cup V(\bar{F_1})\subseteq S_1\cup S$. Therefore, $d_{G_P}(u)+d_{G_P}(v)\le \vert S_1\vert+\vert S\vert=2k-1$, which contradicts to $d_{G_P}(u)+d_{G_P}(v)\ge2(k+\lfloor\frac{m}{2}\rfloor)$.

Therefore $F\cap \bar{F_1}\neq\emptyset$ if $(\alpha)$ holds;  $|S(F,\bar{F_1})|\geq k$ if $(\beta)$ holds. By Lemma 2.2, we have $F_1\cap \bar{F}_{G_P}=\emptyset$,  and thus $\bar{F}_{G_P}\subseteq S_1$. Take $u\in \bar{F}_{G_P}$, since $N_G(u)\subseteq S\cup V(P)\cup V(\bar{F}_{G_P})$ and $u$ has at most $\lceil\frac{|P|}{2}\rceil$ neighbors in $V(P)$, we have   $d_{G[\bar{F}_{G_P}]}(u)\geq k+m-\lceil\frac{|P|}{2}\rceil-|S| \geq\lfloor\frac{m}{2}\rfloor$. Suppose $d_{G[\bar{F}_{G_P}]}(u)>0$. Then there exists an edge $uv\in E(G[\bar{F}_{G_P}])$. If $(\alpha)$ holds, then $N_{G_P}(u)\cap N_{G_P}(v)=\emptyset$ is obvious. If $(\beta)$ holds, then since $V(C)\subseteq V((F\cap F_1)\cup S(F,F_1))\subseteq F\cup S$, we have $V(\bar{F}_{G_P})\cap V(C)=\emptyset$ and $\{u,v\}\subseteq V(G)\setminus V(C)$. By the definition of $\mathcal{K}_k^b(m)$, we also have $N_{G_P}(u)\cap N_{G_P}(v)=\emptyset$. Moreover, $N_{G_P}(u)\cup N_{G_P}(v)\subseteq S\cup V(\bar{F}_{G_P})\subseteq S\cup S_1$. But then $2(k+\lfloor\frac{m}{2}\rfloor)\leq d_{G_P}(u)+d_{G_P}(v)\le \vert S\vert+\vert S_1\vert=2k-1$, a contradiction.
Therefore $d_{G[\bar{F}_{G_P}]}(u)=0$. Now, by $|P|\leq m$, $d_G(u)\geq\delta(G)\geq k+m$ and $N_G(u)\subseteq S\cup V(P)\cup V(\bar{F}_{G_P})$, we conclude that $m=1$ and $u$ is adjacent to all vertices in $S$ and the only vertex $p$ in $P$. Since  $N_{G}(u)\cap N_{G}(p)=\emptyset$ and $N_{G}(u)\cup N_{G}(p)\subseteq S\cup V(P)\cup V(\bar{F}_{G_P})$, we obtain $|N_G(p)|\leq |V(\bar{F}_{G_P})|\leq|S_1|=k-1$, a contradiction.
$\Box$

\section{Main results}

\begin{thm}
For all $(G,C)\in\mathcal{K}_{k^+}^{b}(m)$ and $p_0\in G-V(C)$, there is a path $P\subseteq G-V(C)$ of order $m$ starting from $p_0$, such that $\kappa(G-V(P))\ge k$ holds.
\end{thm}

\noindent{\bf Proof.}
We prove the theorem by induction on the order of the graph at the same time for all $m$. For any graph $G$ and a complete subgraph $C$ of $G$ in a pair $(G,C)\in \mathcal{K}_{k^+}^{b}(m)$, each pair of adjacent vertices $u_1,v_1\in V(G)\setminus V(C)$ satisfies $N_G(u_1)\cap N_G(v_1)=\emptyset$ by the definition of $\mathcal{K}_{k^+}^{b}(m)$. Since $\delta_{G}{(G-C)}\ge k+m$, both $u_1$ and $v_1$ have at least $k+m$ neighbors. Thus $|V(G)|\geq |N_G(u_1)\cup N_G(v_1)|=|N_G(u_1)|+|N_G(v_1)|\geq2(k+m)$. Suppose $|V(G)|=2(k+m)$. Then $V(G)= N_G(u_1)\cup N_G(v_1)$. Let $N_G(u_1)=\{v_1,v_2,\cdots,v_{k+m}\}$ and $N_G(v_1)=\{u_1,u_2,\cdots,u_{k+m}\}$. By the definition of $\mathcal{K}_{k^+}^{b}(m)$, each vertex $u\in\{u_1,u_2,\cdots,u_{k+m}\}\setminus V(C)$ is adjacent to any vertex in $\{v_1,v_2,\cdots,v_{k+m}\}$, and each vertex $v\in\{v_1,v_2,\cdots,v_{k+m}\}\setminus V(C)$ is adjacent to any vertex in $\{u_1,u_2,\cdots,u_{k+m}\}$. Therefore, $G$ is isomorphic to $K_{k+m,k+m}\langle S\rangle$ for some $S\subseteq V(K_{k+m,k+m})$ with $|S|=k$. This implies that the smallest graph $G$ in a pair $(G,C)\in \mathcal{K}_{k^+}^{b}(m)$ is isomorphic to $K_{k+m,k+m}\langle S\rangle$ for some $S\subseteq V(K_{k+m,k+m})$ with $|S|=k$. Since $K_{k+m,k+m}\langle S\rangle$ is $(k+m)$-connected, we have $K_{k+m,k+m}\langle S\rangle-V(P)$ is $k$-connected for any path $P$ of order $m$ in $G[V(K_{k+m,k+m})\setminus S]$. Thus Theorem 3.1 is true for $(G,C)\in\mathcal{K}_{k^+}^{b}(m)$ with the smallest $G$.

Assume Theorem 3.1 is not true for an integer $k$. Let $G$ be a graph of least order such that for certain $C\subseteq G$, $m$, $p_0\in G-V(C)$, the pair $(G,C)\in \mathcal{K}_{k^+}^{b}(m)$ is a counterexample to Theorem 3.1 for $k$.

Since $\kappa{(G)}\ge k+1$ by the definition of $\mathcal{K}_{k^+}^{b}(m)$, we have $\kappa{(G-p_0)}\ge k$. So we can choose a path $P\subseteq G-V(C)$ starting from $p_0$ of maximal order with $|P|\geq1$, such that $G_P=G-V(P)$ is $k$-connected. Since Theorem 3.1 is not true for $(G,C),p_0,m$, we have $1\le\vert P\vert<m$, in particular $m\geq2$. Let $p$ be the other end of $P$. Since $G-C$ is a bipartite graph and $P$ is a path in $G-C$, we have $\vert N_G(p)\cap (G-V(C\cup P))\vert\ge k+m-k-\lceil \frac{|V(P)|}{2}\rceil\geq m-\lceil \frac{m-1}{2}\rceil\ge1$. Therefore $p$ has a neighbor $q\in V(G)-V(C\cup P)$. Since $P_1:=P\cup \{pq\}\subseteq G-V(C)$ is a path starting from $p_0$ of order $|P|+1\leq m$, we conclude $\kappa{(G-V(P_1))}< k$ by the choice of $P$. Hence $\kappa (G_P)=k$. Since each pair of adjacent vertices in $V(G-C)$ has no common neighbors in $G$ and $\delta_{G}{(G-C)}\ge k+m$, we have $|V(G)|\geq2(k+m)$ and  $\vert V(G_P-C)\vert\ge 2(k+m)-|V(P)|-|V(C)|\geq k+m+1$. Thus $G_P$ is not complete. Choose a minimum vertex-cut $S$ of $G_P$ such that $G_P-S$ contains a minimal fragment $F$ of $G_P$ with $F\cap C=\emptyset$. Let $\bar{F}_{G_P}=G_P-(V(F)\cup S)$ be the complementary fragment of $F$ to $S$ in $G_P$. If $\vert F\vert=1$ holds, then for the only vertex $u$ in $F$, we have
$$k=\vert S\vert= d_{G_P}(u)\ge k+m-\lceil \frac{|V(P)|}{2}\rceil\geq k+m-\lceil \frac{m-1}{2}\rceil\ge k+1,$$
which is a contradiction. This implies $\vert F\vert\ge2$, and we have $\kappa{(G_P\langle S\rangle-V(\bar{F}_{G_P}))}\ge k+1$ by Lemma 2.1.

Since $\kappa{(G)}>k$, there is a vertex in $P$ adjacent to some vertices in $F$. Let $p^{'}$ be the vertex nearest to $p$ on $P$ with $N_G(p^{'})\cap F\ne \emptyset$. Choose a vertex $q'\in N_G(p')\cap F$. Define $P^{'}=P[p_0,p^{'}]$, $G_{P^{'}}=G-V(P^{'})$, $R=V(P(p^{'},p])$ and $r=\vert R\vert$.  For any $u\in R$, we have $\vert N_G(u)\cap V(G_{P'})\vert\ge k+m-\lceil \frac{|V(P')|}{2}\rceil\ge k+1$. Furthermore, since $P'\subseteq P$ and $G_P$ is $k$-connected, we obtain $G_{P^{'}}$ is $k$-connected. We have $N_G(R)\cap F=\emptyset$ by the choice of $p^{'}$. Hence $S$ is also a minimum vertex-cut of $G_{P^{'}}$ and $F$ is also a minimal fragment of $G_{P^{'}}$. This implies $\kappa(G_{P^{'}})=k$. Let $\bar{F}_{G_{P'}}$ be the complementary fragment of $F$ to $S$ in $G_{P'}$.
Since $\vert F\vert\ge2$, by applying Lemma 2.1 to $F$ in $G_{P'}$, we obtain that $G_{P^{'}}\langle S\rangle[F\cup S]=G_{P'}\langle S\rangle-V(\bar{F}_{G_{P'}})$ is $(k+1)$-connected. Furthermore, $C\subseteq G_{P^{'}}-V(F)$ and
$$
\delta_{G_{P^{'}}}(G_{P^{'}}-V(C))\ge k+m-\lceil\frac{|V(P)|-r}{2}\rceil\ge k+m-\lceil\frac{m-1-r}{2}\rceil\geq k+m-(m-1-r)=k+r+1.
$$
So we have $(G_{P^{'}},C)\in\mathcal{K}_{k}^b(r+1)$. Moreover, by  $\kappa{(G_{P^{'}}\langle S\rangle[F\cup S])}\ge k+1$, we get $(G_{P^{'}}\langle S\rangle[F\cup S],K(S))\in\mathcal{K}_{k^+}^{b}(r+1)$.
Since $G[F\cup S]$ is a subgraph of $G$ and $\vert F\cup S\vert<\vert V(G)\vert$, by the choice of $G$, we can find a path $Q\subseteq F$ of order $r+1$ starting from  $q'$ so that $G_{P^{'}}\langle S\rangle[F\cup S]-V(Q)$ is $k$-connected. Considering the complementary fragment $\bar{F}_{G_{P'}}$  of $F$ to $S$ in $G_{P^{'}}$ to use Lemma 2.3($\beta$),  we can come to a conclusion that $G_{P^{'}}-V(Q)$ is $k$-connected. Define the path $P_2:=P^{'}\cup \{p'q'\}\cup Q\subseteq G-V(C)$. Then $P_2$ has order $\vert P\vert-r+\vert Q\vert=\vert P\vert+1\le m$. But now $G$ has a path $P_2$ starting from $p_0$ such that $|V(P_2)|=|V(P)|+1$ and $G-V(P_2)=G_{P^{'}}-V(Q)$ is $k$-connected, which contradicts to the choice of $P$. $\Box$

Since the proof of Theorem 3.2 is similar to the one for  Theorem 3.1, we will refer to the proof of Theorem 3.1 in the proof of Theorem 3.2.

\begin{thm}
Every $k$-connected bipartite graph $G$ with $\delta {(G)}\ge k+m$ for positive integers $k, m$ contains a path $P$ of order $m$ such that $G-V(P)$ remains $k$-connected.
\end{thm}

\noindent{\bf Proof.} We consider two cases in the following.

\noindent{\bf Case 1.} $\kappa(G)=k$.

Choose a minimum vertex-cut $S$ and then take a minimal fragment $F$ of $G$ contained in 
$G-S$. Thus $\vert S\vert=k$ by $\kappa(G)=k$. Since  $\delta {(G)}\ge k+m$ and $m$ is positive, we have $|F|\geq2$. By Lemma 2.1, $G\langle S\rangle-V(\bar{F})$ is $(k+1)$-connected and   $(G\langle S\rangle-V(\bar{F}), K(S))\in\mathcal{K}_{k^+}^{b}(m)$. By Theorem 3.1, we can find a path $P\subseteq F$ of order $m$ such that $G\langle S\rangle-V(\bar{F}\cup P)$ is $k$-connected. Then $\kappa(G-V(P))\ge k$ holds by Lemma 2.3($\alpha$).

\noindent{\bf Case 2.} $\kappa(G)\ge k+1$.

By contradiction, assume the theorem is false. As in the proof of Theorem 3.1, choose a longest path $P$. Then $1\le\vert P\vert< m$ and $G_P=G-V(P)$ is $k$-connected. Moreover, $\kappa(G_P)=k$. Let $p_0$ and $p$ be the origin and terminus, respectively, of $P$. Choose a minimum vertex-cut $S$ and then take a minimal fragment $F$ of $G_p$ contained in
$G_p-S$. By $\kappa(G_P)=k$, we have $|S|=k$. If $\vert F\vert=1$, then for $u\in F$, we have $d_{G_P}(u)= k$. But $d_{G_P}(u)\ge k+m-\lceil \frac{|P|}{2}\rceil\ge k+1$, a contradiction. Thus $\vert F\vert\ge 2$ and  $\kappa(G_P\langle S\rangle-V(\bar{F}_{G_P}))\ge k+1$ by Lemma 2.1. Let the definitions of $p^{'}, q', P^{'}, R, r$ and $G_{P^{'}}$ be the same as those in the proof of  Theorem 3.1. We have $\kappa(G_{P^{'}})= k$ and $\delta(G_{P^{'}})\ge k+r+1$. Since $N_G(R)\cap F=\emptyset$, $F$ is also a minimal fragment of $G_{P^{'}}$ to $S$.  We obtain $(G_{P^{'}}\langle S\rangle[F\cup S],K(S))\in \mathcal{K}_{k^+}^{b}(r+1)$ by $\kappa{(G_{P^{'}}\langle S\rangle[F\cup S]})\ge k+1$ and
$\delta_{G_{P^{'}}\langle S\rangle[F\cup S]}(G_{P^{'}}\langle S\rangle[F\cup S]-K(S))\ge k+r+1$.
By almost the same argument used in the proof of Theorem 3.1, there is a path $Q$ of order $r+1$ starting from $q'$ in $F$ so that $G_{P^{'}}\langle S\rangle[(F\cup S)-V(Q)]$ is still $k$-connected. By applying Lemma 2.3($\alpha$) to the fragment $\bar{F}_{G_{P'}}$, we have $\kappa(G_{P^{'}}-V(Q))\ge k$. Note that $G_{P^{'}}-V(Q)=G-(V(P^{'})\cup V(Q))$. Then we find a path $P_2:=P^{'}\cup \{p'q'\}\cup  Q$ starting from $p_0$ such that  $\kappa(G-V(P_2))\ge k$ and $|V(P_2)|=|V(P)|+1$, which contradicts the choice of $P$.
$\Box$

\section{Concluding Remarks}

We believe that the condition $\delta {(G)}\ge k+m$ in Theorem 3.2 is not best possible. Maybe $\delta {(G)}\ge k+\lceil\frac{m}{2}\rceil$ is enough. However, it seems difficult to improve the bound $\delta {(G)}\ge k+m$ in Theorem 3.2 by using the present method in this paper. For general tree $T=T[X,Y]$, we propose the following conjecture.

\begin{conj} For every positive integer $k$ and every finite tree $T$ with bipartition $X$ and $Y$ (denote $t=max\{|X|,|Y|\}$),
every $k$-connected bipartite graph $G$ with $\delta(G)\geq k+t$ contains a subgraph $T'\cong T$ such that $\kappa(G-V(T'))\geq k$.
\end{conj}

Since $K_{k+t-1,k+t-1}$ contains no $T'\cong T$ such that $\kappa(G-V(T'))\geq k$, the bound $\delta(G)\geq k+t$ would be best possible when Conjecture 3 were true.

{\flushleft\textbf{Acknowledgements}}

The authors would like to thank the editor and the anonymous reviewers  for their valuable and kind suggestions which greatly improved the original manuscript.

\end{document}